\documentclass[12pt,leqno]{article}
\usepackage{graphicx}
\usepackage{amsmath}
\usepackage{amsfonts}
\usepackage{amssymb}
\usepackage{amscd}
\usepackage{amsthm}

\theoremstyle{plain}
\newtheorem{theorem}{Theorem}

\newtheorem{lemma}[theorem]{Lemma}
\newtheorem{proposition}[theorem]{Proposition}

\theoremstyle{definition}

\newtheorem{example}{Example}

\date{}

\def\p{\pi}

\def\a{\alpha}
\def\aa{\eta}
\def\b{\beta}
\def\bb{\gamma}

\def\la{\bar{\alpha}}
\def\lb{\bar{\beta}}
\def\lbb{\bar{\gamma}}

\def\tr{\mathbf{t}}
\def\t{t}
\def\g{\xi}

\def\e{\varepsilon}

\def\ps{\varphi}
\def\f{\psi}

\def\r{\rho}
\def\O{\Omega}

\def\S{\mathbf{S}}

\def\MCG{\textup{MCG}}
\def\SL{\textup{SL}}
\def\PMCG{\textup{PMCG}}

\def\E{\mathcal{E}}
\def\K{\mathcal{K}}

\def\Z{\mathbb{Z}}
\def\gcd{\textup{gcd}}

\begin{document}

\title{Representations of $(1,1)$-knots}

\author{Alessia Cattabriga \and Michele Mulazzani}

\maketitle

\begin{abstract}
We present two different representations of $(1,1)$-knots and
study some connections between them. The first representation is
algebraic: every $(1,1)$-knot is represented by an element of the
pure mapping class group of the twice punctured torus
$\PMCG_2(T)$. Moreover, there is a surjective map from the kernel
of the natural homomorphism $\Omega:\PMCG_2(T)\to \MCG(T)\cong
\SL(2,\Z)$, which is a free group of rank two, to the class of all
$(1,1)$-knots in a fixed lens space. The second representation is
parametric: every $(1,1)$-knot can be represented by a 4-tuple
$(a,b,c,r)$ of integer parameters, such that $a,b,c\ge 0$ and
$r\in\Z_{2a+b+c}$. The strict connection of this representation
with the class of Dunwoody manifolds is illustrated. The above
representations are explicitly obtained in some interesting cases,
including two-bridge knots and torus knots.\end{abstract}

\noindent {\it Mathematics Subject
Classification 2000:} Primary 57M25, 57N10; Secondary  20F38, 57M12.\\
{\it Keywords:} $(1,1)$-knots, mapping class groups, cyclic
branched coverings, Dunwoody manifolds, Heegaard diagrams, torus
knots.

\section{Introduction and preliminaries} \label{intro}

\begin{figure}[ht]
\begin{center}
\includegraphics*[totalheight=3cm]{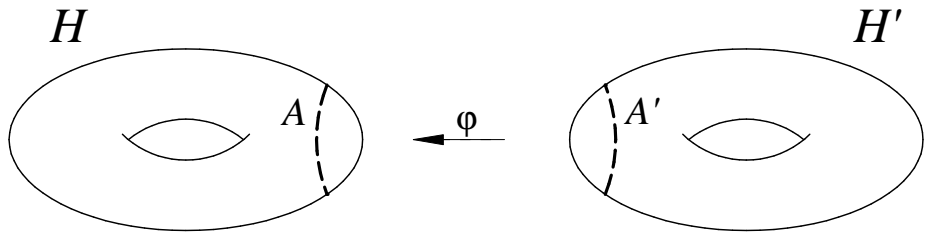}
\end{center}
\caption{A $(1,1)$-knot decomposition.} \label{Fig. 1}
\end{figure}

A knot $K$ in a closed, connected, orientable 3-manifold $N^3$ is
called a $(1,1)$-{\it knot\/} if there exists a Heegaard splitting
of genus one \hbox{$(N^3,K)=(H,A)\cup_{\ps}(H',A'),$} where $H$
and $H'$ are solid tori, $A\subset H$ and $A'\subset H'$ are
properly embedded trivial arcs\footnote{This means that there
exists a disk $D\subset H$ (resp. $D'\subset H'$) with $A\cap
D=A\cap\partial D=A$ and $\partial D-A\subset\partial H$ (resp.
$A'\cap D'=A'\cap\partial D'=A'$ and $\partial
D'-A'\subset\partial H'$).}, and $\ps:(\partial H',\partial
A')\to(\partial H,\partial A)$ is an attaching homeomorphism (see
Figure \ref{Fig. 1}). Obviously, $N^3$ turns out to be a lens
space $L(p,q)$, including $\S^3=L(1,0)$ and
$\S^1\times\S^2=L(0,1)$.

It is well known that the family of $(1,1)$-knots contains all
torus knots and all two-bridge knots in $\S^3$. Several
topological properties of $(1,1)$-knots have recently been
investigated in many papers (see references in \cite{CM2}).

Two knots $K\subset N^3$ and $\bar K\subset\bar N^3$ are said to
be {\it equivalent} if there exists a homeomorphism $f:N^3\to\bar
N^3$ such that $f(K)=\bar K$.

An $n$-fold cyclic covering $M^3$ of a 3-manifold $N^3$, branched
over a knot $K\subset N^3$, is called {\it strongly-cyclic\/} if
the branching index of $K$ is $n$. This means that the fiber in
$M^3$ of each point of $K$ consists of a single point. Observe
that a cyclic branched covering of a knot $K$ in $\S^3$ is always
strongly-cyclic and is uniquely determined, up to equivalence,
since $H_1(\S^3-K)\cong\Z$. Obviously, this property is no longer
true for a knot in a more general 3-manifold.

The necessary and sufficient conditions for the existence and
uniqueness of strongly-cyclic branched coverings of $(1,1)$-knots
have been obtained in \cite{CM1}.

In this paper we present two different representations of
$(1,1)$-knots, as developed in \cite{CM1,CM2,CM3}, and provide new
results.

In Section \ref{algebraic} we show an algebraic representation,
introduced in \cite{CM1,CM2}, through the pure mapping class group
of the twice punctured torus $\PMCG_2(T)$, where $T=\partial H$.
Moreover, we give the proof that the kernel of the natural
homomorphism $\Omega:\PMCG_2(T)\to\MCG(T)\cong\SL(2,\Z)$, is a
free group of rank two. Since there is a surjective map from
$\ker\Omega$ to the class of all $(1,1)$-knots in a fixed lens
space, every $(1,1)$-knot can be represented by an element of
$\ker\Omega$, whose standard generators $\tau_m$ and $\tau_l$ have
a nice topological meaning. A characterization of the subgroup
$\E$ of $\PMCG_2(T)$, consisting of the (isotopy class of)
homomorphisms which extend to the handlebody $H$, fixing $A$, is
also given. The group $\E$ contains elements all producing the
trivial knot in $\S^1\times\S^2$, so its determination appears to
be important in order to produce a ''more injective''
representation.

In Section \ref{parametric} we describe the parametric
representation by 4-tuples of integers, introduced in \cite{CM3}.
This parametrization has a strict connection with the class of
Dunwoody manifolds.

A direct connection between the two representations has been
established in \cite{CM3} for the interesting case of torus knots.
Using this result, an explicit parametrization for a large class
of torus knots is obtained (see Proposition \ref{proptorus}) and a
table with the parametrization for other torus knots is provided
in the Appendix.

\section{Algebraic representation of $\mathbf{(1,1)}$-knots} \label{algebraic}

The mapping class group of a torus $T$ (i.e. the group of the
isotopy class of orientation-preserving homeomorphism of $T$) is
indicated by $\MCG(T)$. Moreover, $\MCG_2(T)$ denotes the mapping
class group of the twice punctured torus, being $P_1$ and $P_2$
two fixed punctures.

Now, let $K\subset L(p,q)$ be a $(1,1)$-knot with
$(1,1)$-decomposition $(L(p,q),K)=(H,A)\cup_{\ps}(H',A')$ and let
$\mu :(H,A)\to (H',A')$ be a fixed orientation-reversing
homeomorphism, then $\f=\ps\mu_{|\partial H}$ is an
orientation-preserving homeomorphism of $(\partial H,\partial
A)=(T,\{P_1,P_2\})$. Moreover, since two isotopic attaching
homeomorphisms produce equivalent $(1,1)$-knots, we have a natural
surjective map $$\Theta:\ \f\in \MCG_2(T)\mapsto K_{\f}\in \K$$
from $\MCG_2(T)$ to the set $\K$ of all $(1,1)$-knots.

In the following, if $\delta$ is a simple closed curve in $T$,
then $\t_{\delta}$ denotes the right-hand Dehn twist around
$\delta$.

\begin{figure}[ht]
\begin{center}
\includegraphics*[totalheight=4cm]{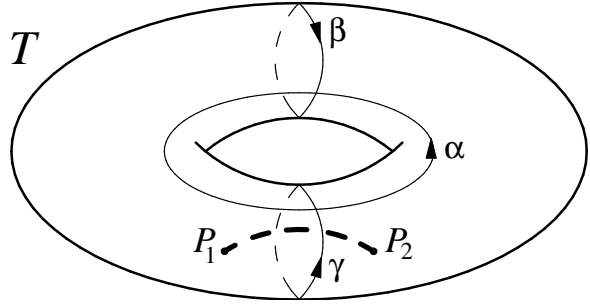}
\end{center}
\caption{Generators of $\PMCG_2(T)$.} \label{Fig. 0}
\end{figure}

Let $\a,\b,\bb$ be the curves depicted in Figure \ref{Fig. 0}.
Then $\MCG_2(T)$ is generated by $\t_{\a},\t_{\b},\t_ {\bb}$, that
fix the punctures, and a $\pi$-radians rotation $\r$ of $T$, that
exchanges the punctures. Observe that $\r$ commutes with the other
generators.

It is easy to see that $\r$ can be extended to a homeomorphism of
the pair $(H,A)$, so $K_{\f}$ and $K_{\f\r}$ are equivalent knots,
for each $\f\in\MCG_2(T)$. Therefore, we can restrict our
attention to the subgroup $\PMCG_2(T)$ of $\MCG_2(T)$, called the
pure mapping class group of the twice punctured torus, consisting
of the elements of $\MCG_2(T)$ fixing the punctures.

The restriction $\Theta'$ of $\Theta$ to $\PMCG_2(T)$ is still
surjective, so every $(1,1)$-knot can be represented by elements
belonging to $\PMCG_2(T)$.

Consider the morphism $\O:\PMCG_2(T)\to\SL(2,\Z)$, obtained as the
composition of the natural epimorphism from $\PMCG_2(T)$ to
$\MCG(T)$ with the isomorphism between $\MCG(T)$ and $\SL(2,\Z)$,
relative to the ordered base $(\b,\a)$ of $H_1(T)$. In terms of
the generators of $\PMCG_2(T)$, $\O$ is given by:
$$\O(t_{\a})= \left(\begin{matrix}1&0\\1&1\end{matrix}\right),\ \
\O(t_{\b})=\O( t_{\bb})=
\left(\begin{matrix}1&-1\\0&1\end{matrix}\right).$$ With the above
notations, if
$\O(\f)=\left(\begin{matrix}q&s\\p&r\end{matrix}\right)$, then
$K_{\f}$ is a $(1,1)$-knot in the lens space $L(\vert p\vert,\vert
q\vert)$ (see \cite[p. 186]{BZ}).

Now we list some examples of $(1,1)$-knots given by this
representation.

\begin{example} \label{example1}
\begin{itemize}
\item[a)] If either $\f=\f_{0,1}=1$ or $\f=t_{\b}$ or
$\f=t_{\bb}$, then $K_{\f}$ is the trivial knot in
$\S^1\times\S^2$. \item[b)] If either $\f=t_{\a}$ or
$\f=\f_{1,0}=\t_{\b}\t_{\a}t_{\b}$, then $K_{\f}$ is the trivial
knot in $\S^3$. \item[c)] Let $p,q$ be integers such that $0<q<p$
and $\gcd(p,q)=1$. If \hbox{$\frac {p}{q}  = a_1+ \frac{1}{a_2
+\frac{1} {a_3 + \cdots + \frac{1}{a_m}}}$,} then the trivial knot
in the lens space $L(p,q)$ is represented by
$$\f_{p,q}=\begin{cases}
t_{\a}^{a_1}t_{\b}^{-a_2}\cdots t_{\a}^{a_m}& \text{ if $m$ is odd}\\
t_{\a}^{a_1}t_{\b}^{-a_2}\cdots t_{\b}^{-a_m}t_{\b}t_{\a}t_{\b} &
\text{if $m$ is even}\end{cases}.$$ \item[d)] If
$\f=t_{\a}t_{\b}t_{\a}t_{\a}t_{\gamma}t_{\a}$, then $K_{\f}$ is
the core knot $\S^1\times\{P\}\subset\S^1\times\S^2$, where $P$ is
any point of $\S^2$.
\end{itemize}
\end{example}

The representation $\Theta'$ is not at all injective and, in
general, there are infinitely many elements of $\PMCG_2(T)$
producing the same $(1,1)$-knot. For example, given $\f\in
\PMCG_2(T)$, all the elements $\f\t_{\b}^c$ produce equivalent
$(1,1)$-knots, for each $c\in\Z$. So a natural question arises: is
it possible to decide if two elements in $\PMCG_2(T)$ represent
the same $(1,1)$-knot? Answering this question seems to be rather
hard.

A first step in this direction is given by the following result.

\begin{theorem} \label{decomposition} {\em \cite{CM2}}
Let $K$ be a $(1,1)$-knot in $L(p,q)$, then there exist
$\f',\f''\in\ker\O$ such that $K=K_{\f}$, with
$\f=\f'\f_{p,q}=\f_{p,q}\f''$, where $\f_{p,q}$ is the map defined
in Example \ref{example1}, only depending on $p$ and $q$.
\end{theorem}

As a consequence, for each lens space $L(p,q)$ we get a surjective
map $$\Theta_{p,q}:\ker\O\to\K_{p,q},$$ where $\K_{p,q}$ is the
set of all $(1,1)$-knots in $L(p,q)$. Moreover, $\ker\O$ has a
very simple structure, as shown in the following result, which is
presented without proof in \cite{CM2}.

\begin{figure}[h]
\begin{center}
\includegraphics*[totalheight=6.5cm]{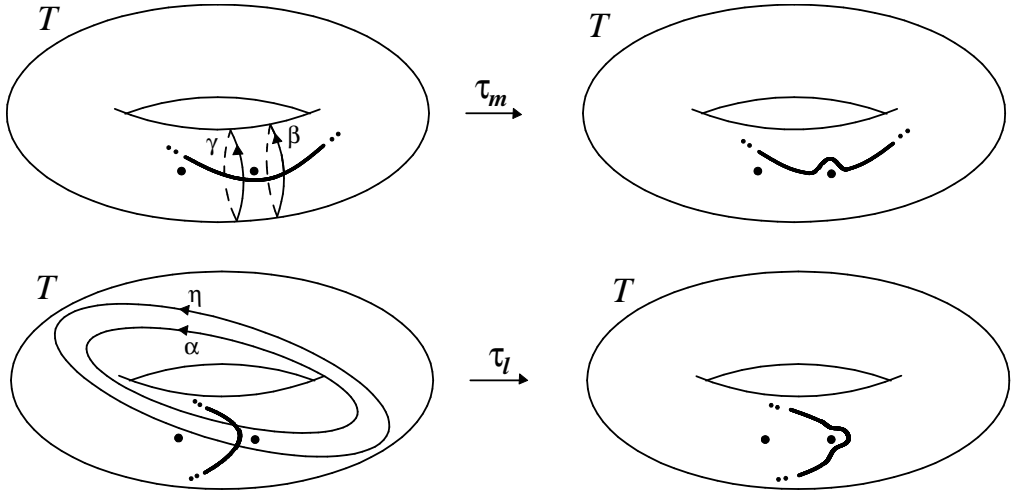}
\end{center}
\caption{Action of $\tau_m$ and $\tau_l$.} \label{Fig. 3}
\end{figure}

\begin{theorem} \label{nucleo} The group $\ker\O$ is freely generated
by $\tau_m=t_{\b}t_{\gamma}^{-1}$ and $\tau_l=t_{\aa}t_{\a}^{-1}$,
where  $t_{\aa}$ is the right-hand Dehn twist around the curve
$\aa$ depicted in Figure \ref{Fig. 3}, and $t_{\aa}=\tau_m^{-1}
t_{\a}\tau_m$.
\end{theorem}

\begin{figure}[h]
\begin{center}
\includegraphics*[totalheight=3.8cm]{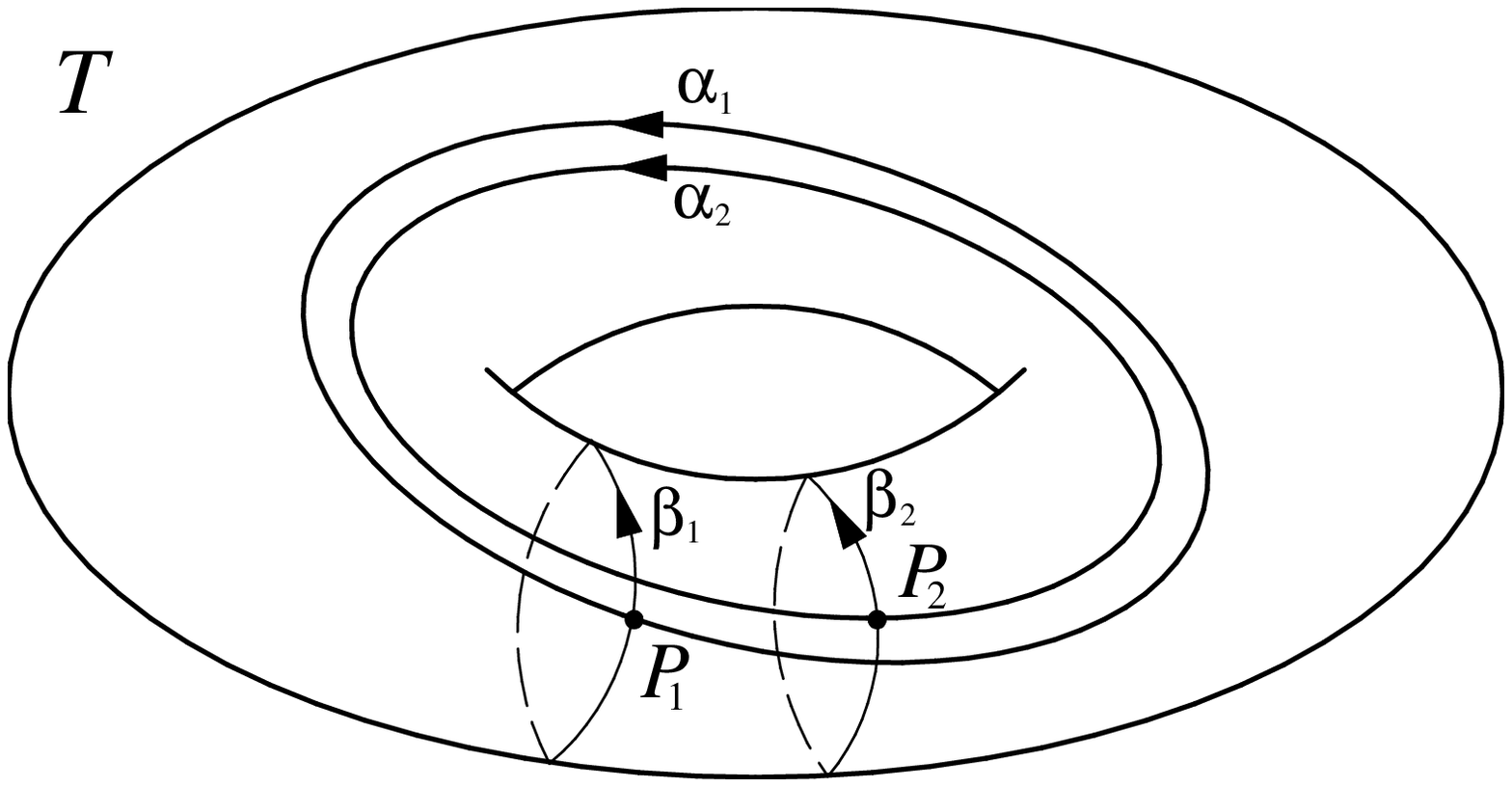}
\end{center}
\caption{} \label{ker}
\end{figure}

\begin{proof}
Let $F_{2}=(T\times T)-\Delta$, where $\Delta=\{(x,x)\,|x\in T\}$
denotes the diagonal, and let $\mathcal{H}(T)$ be the group of
orientation-preserving automorphisms of the torus. Moreover, let
$\mathcal{H}_2$ be  the subgroup of $\mathcal{H}(T)$ consisting of
the elements  pointwise fixing the punctures. By \cite[Th.
1]{Bi2}, the evaluation map $e:\mathcal{H}(T)\rightarrow F_{2}$ is
a fibering with fiber $\mathcal{H}_2$ that induces the exact
sequence on the homotopy groups
$$\cdots\rightarrow\pi_1(\mathcal{H}(T),id)
\stackrel{e_{\#}}{\rightarrow}\pi_1(F_2,(P_1,P_2))\stackrel{d_{\#}}{\rightarrow}\pi_0(\mathcal{H}_2,id)
\stackrel{i_{\#}}{\rightarrow}\pi_0(\mathcal{H}(T),id)\rightarrow
1$$ where $i_{\#}$ denotes the homomorphism induced  by the
inclusion. Since $\pi_0(\mathcal{H}_2,id)=\PMCG_2(T)$ and
$\pi_0(\mathcal{H}(T),id)=\MCG(T)$, we have
$$\ker\O\cong\ker i_{\#}=\textup{im}\,d_{\#}
\cong\pi_1(F_2,(P_1,P_2))/\ker d_{\#}.$$ Moreover, from \cite[Th.
5]{Bi1} we have
$$\pi_1(F_2,(P_1,P_2))=\langle\,\bar{\a_1},\bar{\a_2},\bar{\b_1},\bar{\b_2}\,|\,
1=[\bar{\a_1},\bar{\a_2}]=[\bar{\b_1},\bar{\b_2}]
=[\bar{\a_1},\bar{\b_j}]=[\bar{\b_1},\bar{\a_j}],\,
j=1,2\,\rangle,$$ where $\bar{\a_1}=(\a_1,\a_2)$,
$\bar{\b_1}=(\b_1,\b_2)$, $\bar{\a_2}=(P_1,\a_2)$,
$\bar{\b_2}=(P_1,\b_2)$ where $\a_i$ and $\b_i$ are the loops
depicted in  Figure \ref{ker} and  $P_1$  denotes the constant
loop based on the point $P_1$.  From \cite[Cor. 1.3]{Bi2}, $\ker
d_{\#}$ is freely generated by $\bar{\a_1}$ and $\bar{\b_1}$. So
$\ker\O$ is the free group generated by $d_{\#}(\bar{\a_2})$ and
$d_{\#}(\bar{\b_2})$, which are respectively $\tau_l$ and
$\tau_m$.
\end{proof}

The standard generators $\tau_m$ and $\tau_l$ of $\ker\O$ have a
concrete topological meaning: the effect of $\tau_m$ and $\tau_l$
is to slide one puncture (say $P_2$) respectively along a meridian
and along a longitude of the torus (see Figure~\ref{Fig. 3}).

Since every two-bridge knot admits a Conway presentation with an
even number of even parameters (see \cite[Exercise 2.1.14]{Ka}),
the following result gives a representation for all two-bridge
knots in $\S^3$. An analogous result for torus knots will be given
in Section \ref{torus}.

\begin{proposition} \label{twobridge} {\em \cite{CM2}}
The two-bridge knot having Conway parameters $[2a_1,2b_1,\ldots
,2a_n,2b_n]$ is the $(1,1)$-knot $K_{\f}$ with:
$$\f=\t_{\b}\t_{\a}t_{\b}\tau_m^{-b_n}t_{\e}^{a_n}\cdots
\tau_m^{-b_1}t_{\e}^{a_1},$$ where
$t_{\e}=\tau_l^{-1}\tau_m\tau_l\tau_m^{-1}$.
\end{proposition}

Observe that the representations $\Theta_{p,q}$ are also not
injective, since $K_{\f}$ and $K_{\f\tau_m^c}$ are equivalent
knots, for all $c\in\Z$.

Another way to obtain a ``more injective'' representation seems to
be the characterization of the subgroup $\E$ of $\PMCG_2(T)$,
consisting of the isotopy classes of the homeomorphisms admitting
an extension to an homeomorphism of $H$ which fixes $A$. For each
$\e\in \E$, the knot $K_{\e}$ is the trivial knot in
$\S^1\times\S^2$. Moreover, $\f$ and $\f\e$ produce equivalent
$(1,1)$-knots, for every $\f\in\PMCG_2(T)$ and $\e\in\E$.
Therefore, there exists an induced surjective map
$$\Theta'':\PMCG_2(T)/\E\to \K,$$
where $\PMCG_2(T)/\E$ is the set of left cosets of $\E$ in
$\PMCG_2(T)$.

The following proposition gives a characterization of the elements
of $\E$ in terms of their action on the fundamental groups of
$T-\{P_1,P_2\}$ and $H-A$. Let $*\in T$ be a base point of
$T-\{P_1,P_2\}$. We define the loops $\la=\g\cdot\a\cdot\g^{-1}$,
$\lb=\g_1\cdot\b\cdot\g_1^{-1}$ and
$\lbb=\g_2\cdot\bb\cdot\g_2^{-1}$, where $\g,\g_1,\g_2$ are paths
connecting $*$ to $\a$, $\b$ and $\bb$ respectively. Obviously,
$\pi_1(T-\{P_1,P_2\},*)$ is freely generated by the set
$\{\la,\lb,\lbb\}$ and $\pi_1(H-A,*)$ is freely generated by the
set $\{\la,\lbb\}$.
\begin{proposition} \label{extending} Let $\f\in \PMCG_2(T)$, then
$\f$ belongs to $\E$ if and only if $i_\#(\f_\#(\lb))=1$, where
$i_\#:\p_1(T-\{P_1,P_2\},*)\to \p_1(H-A,*)$ is the homomorphism
induced by the inclusion $i:T-\{P_1,P_2\}\to H-A$.
\end{proposition}
\begin{proof} $\Rightarrow$) Trivial.
$\Leftarrow$) By the proof of \cite[Theorem 10.1]{Gr}, $\f$
extends to a homeomorphism $\widetilde{\f}$ of $H$. Moreover,
$\f(\b)$ bounds a disk $D$ such that $D\cap A=D\cap
\widetilde{\f}(A)=\emptyset$, and the cutting of $H$ along $D$
produces a 3-ball. Therefore, up to isotopy we can suppose that
$\widetilde{\f}(A)=A$.
\end{proof}
It is easy to verify that  $t_{\b},t_{\gamma}$ and
$(t_{\b}t_{\a}t_{\b})^2$ belong to $\E$, while $\t_{\a}$ does not,
but the problem of finding a (possibly finite) presentation for
$\E$ is still open.

\section{Parametric representation of $\mathbf{(1,1)}$-knots} \label{parametric}

As proved in \cite{CM3}, a $(1,1)$-knot $K_{\f}$ is completely
determined by the curve $\f(\beta)$ on $T-\{P_1,P_2\}$. Moreover,
in the open Heegaard diagram obtained by cutting $T$ along
$\beta$, the curve  $\f(\beta)$ is, up to Singer moves \cite{Si}
fixing the set $\{P_1,P_2\}$, one of the three types depicted in
Figure \ref{Fig4} (see proof of \cite[Theorem 3]{CM3}). In all the
cases the circles $C'$ and $C''$ represent the curve $\beta$.

In case (1), the parameters $a,b$ and $c$ denote $a,b$ and $c$
parallel arcs respectively, which are $\f(\b)$ after the cutting.
In this case, we have $d=2a+b+c>0$. The parameter $r$ gives the
gluing rule between the circles $C'$ and $C''$. Obviously, $r$ can
be taken mod $d$. The corresponding $(1,1)$-knot is denoted by
$K(a,b,c,r)$.

In case (2), the corresponding $(1,1)$-knot is the trivial knot in
$\S^1\times\S^2$, which is denoted by $K(0,0,0,0)$.

In case (3), the corresponding $(1,1)$-knot is the core knot
$\S^1\times\{P\}\subset\S^1\times\S^2$, which admits no
parametrization, as will be explained in the following.

In this way we obtain a parametrization of $(1,1)$-knots by a
4-tuple of integers $(a,b,c,r)$, with $a,b,c\ge 0$, and either
$r\in\Z_{d}$, when $d>0$, or $r=0$, when $d=0$.

An interesting property of this parametrization is its connection
with Dunwoody manifolds, which are closed orientable 3-manifolds
introduced in \cite{Du} using a  class of trivalent regular planar
graphs (called {\it Dunwoody diagrams\/}), depending on six
integer parameters $a,b,c,n,r,s$, such that $n>0$, $a,b,c\ge 0$.

\begin{figure}
 \begin{center}
 \includegraphics*[totalheight=6.5cm]{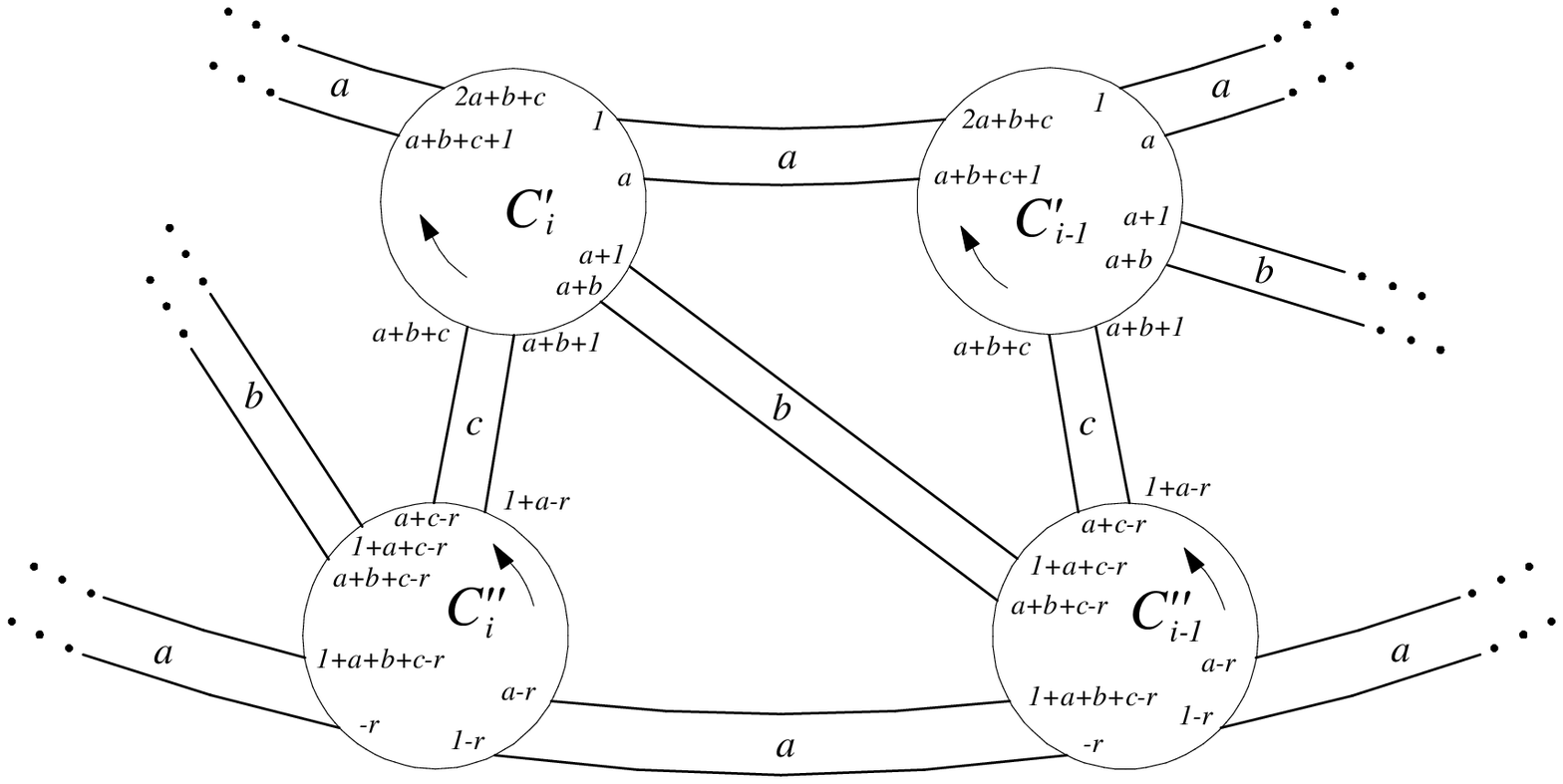}
 \end{center}
 \caption{Heegaard diagram of Dunwoody type.}
 \label{Fig. 13}
\end{figure}

\begin{figure}[h]
\begin{center}
\includegraphics*[totalheight=9cm]{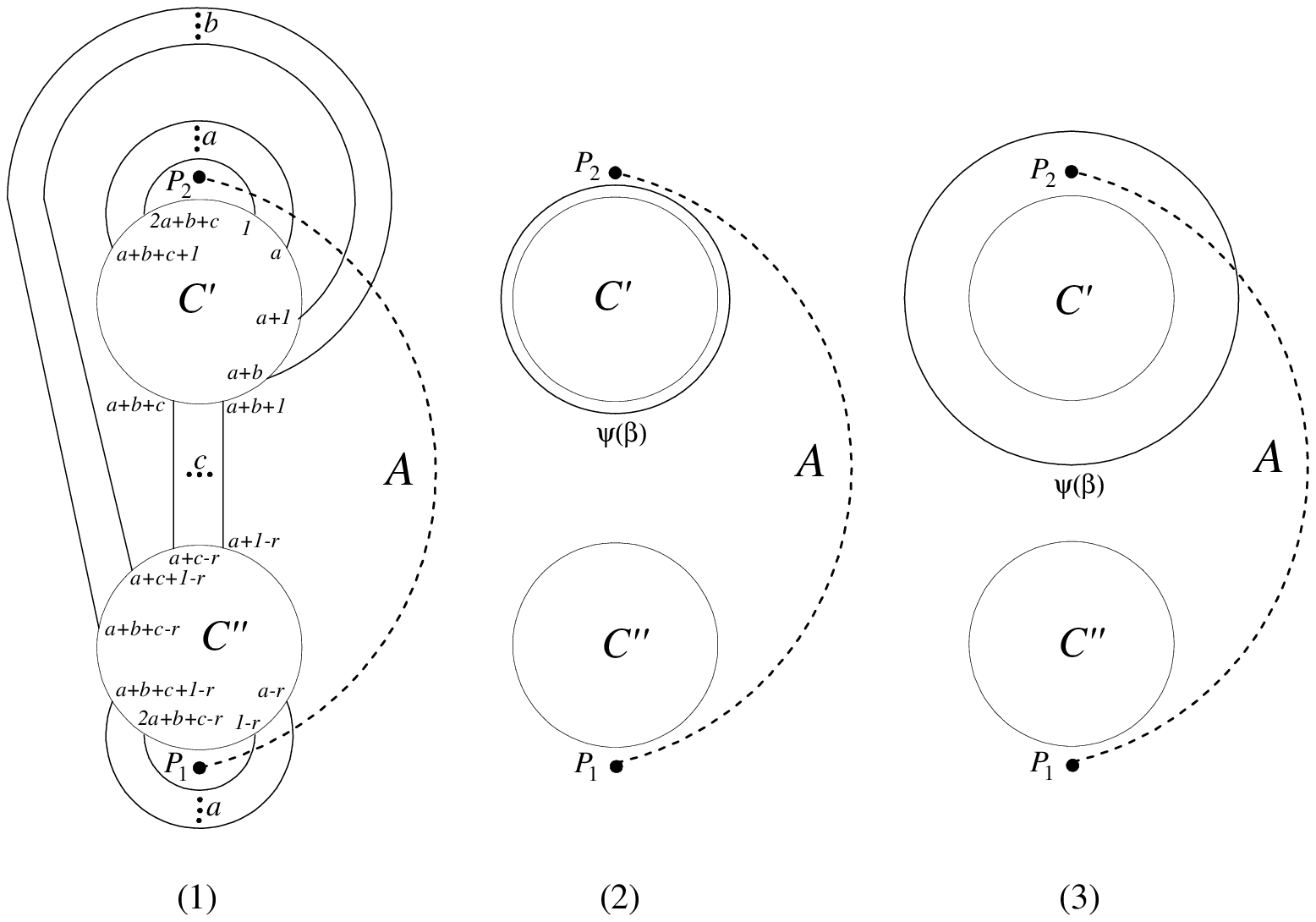}
\end{center}
\caption{} \label{Fig4}
\end{figure}

More precisely, for particular values of the parameters, called
{\it admissible}, a Dunwoody diagram is an (open) Heegaard diagram
of genus $n$ (see Figure \ref{Fig. 13}), which contains $n$
internal circles $C'_1,\ldots,C'_n$, and $n$ external circles
$C''_1,\ldots,C''_n$, each having $d=2a+b+c$ vertices. For every
$i=1,\ldots,n$, the circle $C'_i$ (resp. $C''_i$) is connected to
the circle $C'_{i+1}$ (resp. $C''_{i+1}$) by $a$ parallel arcs, to
the circle $C''_{i}$ by $c$ parallel arcs and to the circle
$C''_{i-1}$ by $b$ parallel arcs (subscripts mod $n$). The cycle
$C'_i$ is glued to the cycle $C''_{i+s}$ (subscripts mod $n$) so
that equally labelled vertices are identified.

Observe that the parameters $r$ and $s$ can be considered mod $d$
and $n$, respectively. Since the identification rule and the
diagram are invariant with respect to an obvious cyclic action of
order $n$, the Dunwoody manifold $D(a,b,c,r,n,s)$ admits a cyclic
symmetry of order $n$.

\begin{theorem} \label{Theorem 3} \begin{itemize} \item[i)] {\em \cite{GM}}
The Dunwoody manifold $D(a,b,c,n,r,s)$ is the $n$-fold
strongly-cyclic covering of the lens space $D(a,b,c,1,r,0)$
(possibly $\S^3$), branched over $K(a,b,c,r)$. \item[ii)] {\em
\cite{CM3}} If $M^3$ is an $n$-fold strongly-cyclic branched
covering of $K(a,b,c,r)$, then there exists $s\in\Z_n$ such that
$M^3$ is homeomorphic to the Dunwoody manifold $D(a,b,c,n,r,s)$.
\end{itemize}
\noindent Therefore, the class of Dunwoody manifolds coincides
with the class of strongly-cyclic branched coverings of
$(1,1)$-knots.
\end{theorem}

The core knot cannot be parameterized as $K(a,b,c,r)$, since it
admits no strongly-cyclic branched coverings (see \cite{CM1}).

Observe that not every 4-tuple of non-negative integers
$(a,b,c,r)$ determines a $(1,1)$-knot $K(a,b,c,r)$, since the
corresponding diagram could fail to be a Heegaard diagram. For
example, the 4-tuples $(a,0,a,a)$, with $a>1$, and $(1,0,c,2)$,
with $c$ even, do not determine any $(1,1)$-knot (see \cite{GM}).

\vbox{
\begin{example} \label{example2}
\begin{itemize}
\item[a)] The trivial knot in $L(p,q)$ (including
$L(1,0)\cong\S^3$) is $K(0,0,p,q)$. \item[b)] The two-bridge knot
of type $(2a+1,2r)$ is $K(a,0,1,r)$ (see \cite{GM}). \item[c)] The
$(1,1)$-knot $K(1,1,1,2)\subset \S^1\times\S^2$ admits three
3-fold strongly-cyclic branched coverings. One of them is the
3-torus $\S^1\times\S^1\times\S^1$, which is homeomorphic to the
Dunwoody manifold $D(1,1,1,3,2,1)$. It is well known that this
manifold cannot be a cyclic branched covering of any knot in
$\S^3$.
\end{itemize}
\end{example}
}

As well as for the algebraic representation, the parametric
representation of a $(1,1)$-knot is not unique, as proved by the
following lemma.

\begin{figure}
\begin{center}
\includegraphics*[totalheight=18cm]{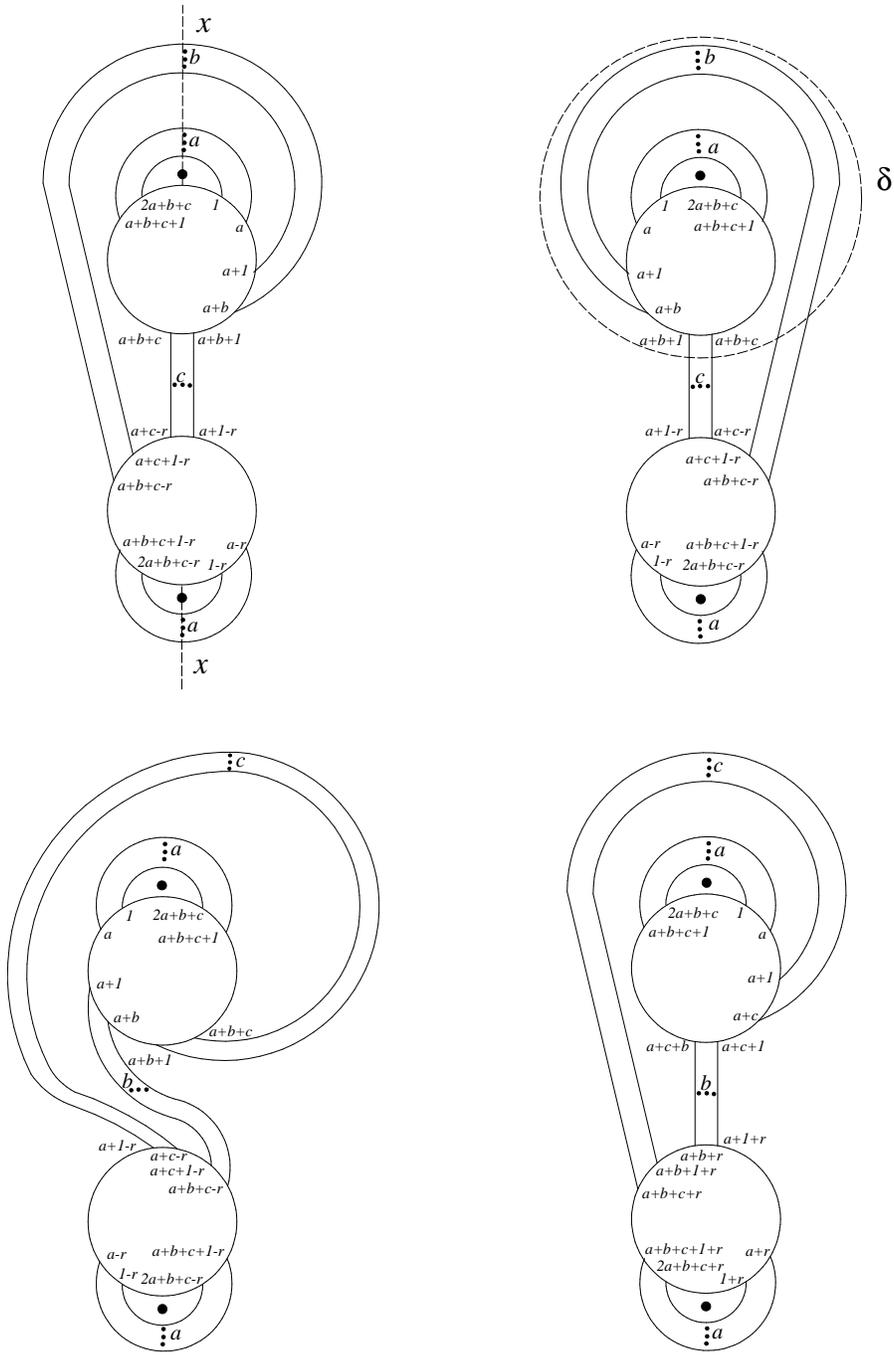}
\end{center}
\caption{From $K(a,b,c,r)$ to $K(a,c,b,-r)$.} \label{Fig.5}
\end{figure}

\begin{lemma}
\begin{itemize}
\item[\textup{a)}]  $K(a,b,c,r)$ and $K(a,c,b,-r)$ are equivalent;
\item[\textup{b)}] $K(a,0,c,r)$ and $K(a,c,0,r)$ are equivalent.
\end{itemize}
\end{lemma}
\begin{proof} a) Looking at Figure \ref{Fig.5}, we pass from the first
diagram, representing $K(a,b,c,r)$, to the second by a reflection
along an axis passing through the punctures (denoted by $x$-$x$ in
the figure). Operating a Singer move of type IIB along $\delta$,
and relabelling the vertices, we obtain $K(a,c,b,-r)$. b) The
application of a Singer move of type IIB along $\delta$ (see
Figure \ref{Fig.5}) on $K(a,b,0,r)$ gives $K(a,0,b,r)$.
\end{proof}

A different parametrization of $(1,1)$-knots, involving four
parameters for the knot and two additional parameters for the
ambient space, can be found in \cite{CK}.

\section{The case of torus knots} \label{torus}
As previously remarked, a very important class of $(1,1)$-knots in
$\S^3$ are torus knots. Without loss of generality, we can
consider torus knots ${\bf t}(k,h)$, with $0<k<h$. The next result
gives the algebraic representation for torus knots. In the
following, $\lfloor x \rfloor$ denotes the integral part of $x$.

\begin{proposition} \label{torusknots} {\em \cite{CM2}} The torus knot ${\bf
t}(k,h)$ is the $(1,1)$-knot $K_{\f}$ with:
\begin{equation} \label{toruseq}
\f=\prod_{j=0}^{h-1}(\tau_l^{-1}\tau_m^{\varepsilon_{h-j}})\t_{\b}\t_{\a}t_{\b},
\end{equation}
where $\varepsilon_{h-j}=\lfloor (j+1)k/h \rfloor-\lfloor (j+2)k/h
\rfloor$, $\tau_m=\t_{\b}\t_{\bb}^{-1}$ and
$\tau_l=\tau_m^{-1}\t_{\a}\tau_m\t_{\a}^{-1}$.
\end{proposition}

Moreover the following proposition tells us how to pass from the
algebraic to the parametric representation of a torus knot.

\begin{proposition} \label{algorithm} {\em \cite{CM3}}
Let $\tr(k,h)\subset\S^3$ be a torus knot and  $\f$ be   as in
\textup{(\ref{toruseq})}. Then $\tr(k,h)=K(a,b,c,r)$, where
$(a,b,c,r)=(a_h,b_h,c_h,r_h)$ is the final step of the following
algorithm, applied for $i=h-j=1,\ldots,h$:
\begin{itemize}
\item[--] $(a_0,b_0,c_0,r_0)=(0,0,1,0)$ and $z_0=0$; \item[--] for
$i=1,\ldots,h$:
$$\left\{\begin{array}{l}a_{i}=a_{i-1}+v\\
b_{i}= r_{i-1}-2w-ud\\
c_{i}= d-b_{i}\\
r_{i}= a_{i-1}+v+w\\
z_{i}= u -\varepsilon_i
\end{array}\right.$$
where: $$ w=\begin{cases}
    a_{i-1}+b_{i-1}+c_{i-1}& \text{if $z_{i-1}<-1-\varepsilon_i$}\\
    a_{i-1}+c_{i-1}& \text{if $z_{i-1}=-1-\varepsilon_i$}\\
     a_{i-1}  & \text{if $z_{i-1}>-1-\varepsilon_i$}
      \end{cases},$$
      $$ v=\begin{cases}
    -(b_{i-1}+c_{i-1})(z_{i-1}+1+\varepsilon_i)-b_{i-1}& \text{if $z_{i-1}<-1-\varepsilon_i$} \\
    0& \text{if $z_{i-1}=-1-\varepsilon_i$}\\
    (b_{i-1}+c_{i-1})(z_{i-1}+1+\varepsilon_i)-c_{i-1}& \text{if $z_{i-1}>-1-\varepsilon_i$}
  \end{cases},$$
 and $u=\lfloor(r_{i-1}-2w)/d\rfloor$, with $d=2a_{i-1}+b_{i-1}+c_{i-1}$.
\end{itemize}
\end{proposition}

Explicit formulae for torus knots of type $\tr(k,ck\pm 1)$ have
been obtained in \cite{AGM,CM3} (see Appendix).

The next result  gives the explicit parametric representation of
another family of torus knots, which contains all the torus knots
with bridge number at most three.

\begin{proposition} \label{proptorus} The torus knot $\tr(sq'+1,(sq'+1)q+s)$,
is
$$K(q',q'(2qq'(s-1)+2q+s-2),1+(s-2)q',2q'^2(s-1)+sq'+1),$$
for every $q,q'>0$ and $s>1$.
\end{proposition}
\begin{proof} From Proposition \ref{torusknots}, we obtain that $\tr(sq'+1,(sq'+1)q+s)$
is  represented by
$\f=((\tau_l^{-q}\tau_m^{-1})^{q'}\tau_l^{-1})^{s-1}(\tau_l^{-q}\tau_m^{-1})^{q'+1}\tau_l^{-1}\t_{\b}\t_{\a}\t_{\b}$.
By \cite[Corollary 8]{CM3}, the application of
$(\tau_l^{-q}\tau_m^{-1})^{q'+1}\tau_l^{-1}$ to $K(0,0,1,0)$ gives
$K(1,q'-1,q'(2q-1),q'+1)$ and $z=0$. Applying $\tau_l^{-1}$, we
get $K(q',q'-1,2+q'(2q-1), q'+1)$ and $z=0$. Then we have to apply
$q'$ times $(\tau_l^{-q}\tau_m^{-1})$. If $q'>1$, applying the
first $(\tau_l^{-q}\tau_m^{-1})$,  we get
$K(q',q'-2,q'(4q-1)+3,3q'+2)$ and $z=0$. So, each time we apply
$(\tau_l^{-q}\tau_m^{-1})$,  under the condition that it is not
the final step, the $a$ and $z$ terms remain unchanged, the $b$
term decreases by one, the $c$ term increases by $2qq'+1$ and the
$r$ term increases by $2q'+1$. So, after $(q'-1)$ steps, we get
$K(q',0,2qq'^2+1,2q'^2)$ and $z=0$. Now applying
$(\tau_l^{-q}\tau_m^{-1})$ for the last time, we get
$K(q',2qq'^2+2qq',1,2q'^2+2q'+1)$ and $z=-1$, which equals the
formula for $s=2$. If $s>2$  we have to apply
$(\tau_l^{-q}\tau_m^{-1})^{q'}\tau_l^{-1}$ to
$K(q',2qq'^2+2qq',1,2q'^2+2q'+1)$, with $z=-1$, another time.
Proceeding as before, we obtain
$K(q',4qq'^2+2qq'+q',q'+1,4q'^2+3q'+1)$ and $z=-1$. So each time
we apply $(\tau_l^{-q}\tau_m^{-1})^{q'}\tau_l^{-1}$, the $a$ and
$z$ terms remain unchanged, the $b$ term increases by $2qq'^2+q'$,
the $c$ term increases by $q'$ and the $r$ term increases by
$2q'^2+q'$. So, after $(s-1)$ steps, we get
$K(q',(2qq'^2+q')(s-2)+2qq'^2+2qq',q'(s-2)+1,(2q'^2+q')(s-2)+2q'^2+2q'+1)$,
as stated.
\end{proof}

The algorithm of Proposition \ref{algorithm} can easily be
implemented. The table in the Appendix is obtained by computer and
contains the parametrization of all torus knots ${\bf t}(k,h)$,
with $k,h\le 25$, not included in the previous cases.

\section{Appendix - $\mathbf{(1,1)}$-parametrization of torus knots} \label{appendix}

\begin{itemize} \item $\tr(k,qk+1)$ is $K(1,k-2,(k-1)(2q-1),k)$,
for all $k>1$ and $q>0$ (see \cite{AGM,CM3}). \item ${\bf
t}(k,qk-1)$ is $K(1,k-2,(k-1)(2q-1)-2,(k-1)(2q-3))$, for all
$k,q>1$ (see \cite{AGM}). \item for all $q,q'>0$ and $s>1$,
$\tr(sq'+1,(sq'+1)q+s)$ is
$$K(q',q'(2qq'(s-1)+2q+s-2),1+(s-2)q',2q'^2(s-1)+sq'+1).$$
\end{itemize}

The following table gives the parametrization of the torus knots
$\tr(k,h)$, with $k,h\le 25$, non included in the previous
formulae.
\bigskip
\begin{center}
\begin{tabular}[tbc]{ccccc}
\begin{tabular}[tbc]{|l|l|}\hline \ &\ \\ \textbf{Knot}&\textbf{Parametrization}\\ \ &\ \\ \hline$\tr(5,8)$ &
$K(2,1,14,11)$ \\\hline $\tr(5,13)$ & $K(2,1,26,11)$\\
\hline$\tr(5,18)$ & $K(2,1,38,11)$\\\hline$\tr(5,23)$ & $K(2,1,50,11)$\\
\hline$\tr(7,11)$ & $K(2,3,24,17)$\\\hline$\tr(7,12)$ & $K(3,1,34,23)$\\\hline$\tr(7,18)$ & $K(2,3,44,17)$\\
\hline$\tr(7,19)$ & $K(3,1,58,23)$\\\hline$\tr(7,25)$ & $K(2,3,64,17)$\\
\hline$\tr(8,11)$ & $K(3,2,33,28)$\\\hline$\tr(8,13)$ & $K(3,41,2,32)$\\\hline$\tr(8,19)$ & $K(3,2,63,28)$\\
\hline$\tr(8,21)$ & $K(3,71,2,32)$\\
\hline$\tr(9,14)$ & $K(2,5,34,23)$\\\hline$\tr(9,16)$ & $K(4,1,62,39)$\\
\hline$\tr(9,23)$ & $K(2,5,62,23)$\\\hline$\tr(9,25)$ & $K(4,1,102,39)$\\
\hline$\tr(10,17)$ & $K(3,4,61,38)$\\

\hline
\end{tabular}&\ &\ &\ &
\begin{tabular}[tbc]{|l|l|}\hline \ &\ \\ \textbf{Knot}&\textbf{Parametrization}\\ \ &\ \\
\hline$\tr(11,14)$ & $K(4,3,60,53)$\\\hline$\tr(11,15)$ & $K(3,5,54,43)$\\
\hline$\tr(11,17)$ & $K(2,7,44,29)$\\\hline$\tr(11,18)$ & $K(3,68,5,53)$\\
\hline$\tr(11,19)$ & $K(4,86,3,59)$\\\hline$\tr(11,20)$ & $K(5,1,98,59)$\\
\hline$\tr(11,25)$ & $K(4,3,116,53)$\\
\hline$\tr(12,17)$ & $K(5,2,87,68)$\\\hline$\tr(12,19)$ & $K(5,99,2,72)$\\
\hline$\tr(13,18)$ & $K(5,98,3,83)$\\
\hline$\tr(13,20)$ & $K(2,9,54,35)$\\\hline$\tr(14,17)$ & $K(5,4,95,86)$\\
\hline$\tr(14,19)$ & $K(3,8,75,58)$\\\hline$\tr(15,19)$ & $K(4,7,96,81)$\\
\hline$\tr(17,20)$ & $K(6,5,138,127)$\\
\hline$\tr(18,23)$ & $K(7,179,4,158)$\\\hline$\tr(18,25)$ & $K(5,163,8,138)$\\
\hline$\tr(19,23)$ & $K(5,9,150,131)$\\
\hline$\tr(19,24)$ & $K(4,11,132,109)$\\
\hline$\tr(20,23)$ & $K(7,6,189,176)$\\
\hline\ &\ \\
\end{tabular}
\end{tabular}
\end{center}

\bigskip

\noindent{\bf Acknowledgements}

\noindent Work performed under the auspices of the G.N.S.A.G.A. of
I.N.d.A.M. (Italy) and the University of Bologna, funds for
selected research topics.

\vspace{15 pt} {ALESSIA CATTABRIGA, Department of Mathematics,
University of Bologna, Piazza di Porta San Donato 5, 40126 Bologna
(Italy). E-mail: cattabri@dm.unibo.it}

\vspace{15 pt} {MICHELE MULAZZANI, Department of Mathematics and
C.I.R.A.M., University of Bologna, Piazza di Porta San Donato 5,
40126 Bologna (Italy). E-mail: mulazza@dm.unibo.it}

\end{document}